\documentclass[11pt,reqno]{article}

\usepackage{amscd,amssymb}
\usepackage[mathscr]{eucal}
\usepackage[all]{xy}
\usepackage{mathrsfs}
\usepackage{xypic}
\usepackage{amsfonts}
\usepackage{amsmath}
\usepackage{amsthm}
\usepackage{latexsym}
\usepackage{eufrak}

\theoremstyle{plain}
\newtheorem{theorem}[equation]{Theorem}

\newtheorem{proposition}[equation]{Proposition}
\newtheorem{lemma}[equation]{Lemma}

\theoremstyle{definition}
\newtheorem{definition}[equation]{Definition}

\newtheorem{remark}[equation]{Remark}

\makeatletter
\renewcommand{\subsection}{\@startsection{subsection}{2}{0pt}{-3ex
plus -1ex minus -0.2ex}{-2mm plus -0pt minus
-2pt}{\normalfont\bfseries}} \makeatother

\numberwithin{equation}{subsection}
\newcommand{\op}{\operatorname}
\def\mto{\longmapsto}

\def\l{\label}

\newcommand{\beq}{\begin{equation}\label}
\newcommand{\eeq}{\end{equation}}

\newcommand{\iso}{{\;\;\stackrel{_\sim}{\longrightarrow}\;\;}}

\newcommand{\vi}{${\sf {(i)}}\;$}
\newcommand{\vii}{${\sf {(ii)}}\;$}
\newcommand{\viii}{${\sf {(iii)}}\;$}

\newcommand{\sset}{\subset}
\DeclareMathOperator{\Ker}{\mathrm{Ker}}
\DeclareMathOperator{\im}{\mathrm{Im}}
\DeclareMathOperator{\Lie}{\mathrm{Lie}}

\def\map{\longrightarrow}
\def\pb{$\bullet\quad$\parbox[t]{142mm}}

\def\ccirc{{{}_{^{\,^\circ}}}}

\newcommand{\too}{\,\,\longrightarrow\,\,}

\newcommand{\inv}{^{-1}}
\newcommand{\ad}{{\mathtt{{ad}}^{\,}}}
\newcommand{\Ad}{{\mathtt{{Ad}}^{\,}}}

\newcommand{\sym}{{\mathtt{{Sym}^{\,}}}}
\newcommand{\hdot}{{\:\raisebox{3pt}{\text{\circle*{1.5}}}}}

\newcommand{\idot}{{\:\raisebox{1pt}{\text{\circle*{1.5}}}}}

\newcommand{\Hom}{{\mathtt{Hom}}}

\newcommand{\gp}{\g_{_P}}
\newcommand{\dg}{_{_{\op{DG}}}}
\newcommand{\Aut}{{\mathtt {Aut}}}

\newcommand{\DT}{{\mathsf{T}}}
\newcommand{\GL}{{\mathsf{GL}}}

\newcommand{\om}{\omega}

\newcommand{\Om}{\Omega}
\newcommand{\wh}{\widehat}

\newcommand{\mc}{{\mathsf{MC}}}
\newcommand{\M}{{\mathscr{M}}}
\newcommand{\OO}{{\mathbb{O}}}
\newcommand{\Si}{\Sigma}

\def\be{\beta}
\def\l{\label}

\def\C{{\mathbb{C}}}

\def\Z{{\mathbb{Z}}}

\def\byy{{\bar{x}}}

\def\G{{\mathfrak{G}}}

\def\g{{\mathfrak{g}}}

\begin{document}
\setlength{\parindent}{6mm}
\setlength{\parskip}{3pt plus 5pt minus 0pt}

\centerline{\Large {\bf Hamiltonian reduction and Maurer-Cartan
equations}}

\vskip 6mm
\centerline{\large {\sc Wee Liang Gan and Victor Ginzburg}}

\bigskip
\centerline{\it To Boris Feigin on the occasion of his 50-th Birthday}
\bigskip

\begin{abstract} We show that solving the Maurer-Cartan equations
is, essentially, the same thing as performing the Hamiltonian
reduction construction. In particular, any
differential graded Lie algebra equipped with an even nondegenerate 
invariant bilinear form gives rise to modular stacks with
symplectic structures.
\end{abstract}

\bigskip
This paper is our modest present to Boris Feigin from whom
we learned the power of Homological Algebra.

\section{From Maurer-Cartan equations to Moment maps}
\setcounter{equation}{0}
\subsection{}\label{sec1}
Let $\Bbbk$ denote the field of either  real
or  complex numbers.
Let $\G=\G_+\oplus\G_-$ be a differential $\Z/2\Z$-graded Lie
(super)-algebra (DGLA) over $\Bbbk$.
Write $d: \G_\pm\to\G_\mp$ for the
differential, and $H_\idot(\G,d):=
\Ker d/\im d$ for the corresponding homology space,
which is again  $\Z/2\Z$-graded Lie
super-algebra.

Depending on the problem, the space $\G$
may be either finite or infinite dimensional.
To fix ideas,
we shall assume below that $\G$ has
finite dimension over $\C$;
moreover, it will be assumed that $\G_+=\Lie G_+$ is
the Lie algebra of a complex connected simply-connected 
linear algebraic group $G_+$. These assumptions,
though certainly too restrictive,
will allow us to make the main ideas more clear.
In reality, the algebra $\G$ itself is  typically
 infinite-dimensional while the homology algebra,
$H_\idot(\G,d),$ is typically finite-dimensional.
In such cases, the space $\G$  usually comes equipped with
a natural topology. Various analytic issues (e.g.,
the closedness of the kernel and image of the
differential $d$) that arise in such a topological
framework require a considerable amount
of machinery and are beyond the scope of this
short paper.

\subsection{}\label{sec2}
The algebraic group $G_+$ acts on its Lie algebra  $\G_+$ via the adjoint action.
Fix an  $\Ad G_+$-orbit $\OO\sset \G_+$, and consider the
following (locally closed) subscheme in $\G_-$:
\[\mc(\G,\OO):=\big\{ x\in \G_-\enspace\big|\enspace dx+\frac{1}{2}[x,x] 
\in \OO\big\}\,. \]
The equation $dx+\frac{1}{2}[x,x]=0$ is known as the {\it Maurer-Cartan}
equation, so we call $\mc(\G,\OO)$ the Maurer-Cartan scheme
\footnote{It is actually a DG scheme, cf. \cite{CFK},
 the fact that will be
exploited
later.}
associated to an orbit $\OO$. If $\OO=\{0\}$ is a one-point orbit,
then our scheme  $\mc(\G,\OO)$ reduces to the zero-scheme of
the above-mentioned standard Maurer-Cartan equation.

Given $a\in \G_+$, let $\xi_a$ denote an  affine-linear
algebraic vector field on $\G_-$ whose value
at the point $x\in \G_-$ is  $\xi_a(x):= [a,x] -da$.
The lemma below is well-known, see e.g. \cite[sect. 1.3]{GM}.
\begin{lemma}\l{action}
\vi The map $a \longmapsto \xi_a$ is a Lie algebra  homomorphism.
 
\vii For any orbit $\OO\sset \G_+$ and any
 $a\in \G_+,$ the vector field $\xi_a$ is
tangent to the Maurer-Cartan scheme $\mc(\G,\OO)\sset \G_-$.
\end{lemma}
It follows from part (i) of the Lemma that, exponentiating 
the vector fields $\xi_a,\,a\in\G_+$, one obtains
an action of the group $G_+$ on $\G_-$ by  affine-linear
transformations. This $G_+$-action on $\G_-$ is known as the
{\em gauge action} (it should not be confused
with the ordinary $\Ad G_+$-action on $\G_-$).
We observe also that (by the Lemma) the map  $a \longmapsto \xi_a$ intertwines the
$\Ad G_+$-action on $\G_+$ with  $G_+$-action on vector fields
induced by the gauge action on $\G_-$.
Further, part (ii) of the Lemma implies
 that the  Maurer-Cartan scheme  $\mc(\G,\OO)$ is stable
under the gauge action of $G_+$.

\begin{definition}\l{modular_stack} We write $\M(\G,\OO):=\mc(\G,\OO)/G_+$
for
the stack-quotient of  $\mc(\G,\OO)$ by the  gauge $G_+$-action,
cf. \cite{LMB},\cite{To},
and call $\M(\G,\OO)$ the {\textsl{modular stack}} attached
to the orbit $\OO\sset \G_+$.
\end{definition}

\begin{remark}\label{deg1} Usually, one has a natural $\Z$-grading
$\G=\bigoplus_i 
\,\G_i$,
making $\G$ a DGLA with differential $\G_\idot\to\G_{\bullet+1}$.
We then  put $\G_+=\bigoplus_i\, \G_{2i}$ and $\G_-=\bigoplus_i\, \G_{2i+1}$.
In such a case, one introduces a smaller Lie group
$G_0\sset G_+$ corresponding to the Lie subalgebra $\G_0\sset \G_+$.
This group acts naturally on $\G_2$, and for any $G_0$-orbit $\OO_2\sset
\G_2$ we may form the corresponding $\G_+$-orbit $\OO=\G_+(\OO_2)$.
Define
\[ \mc_1(\G,\OO):=\big\{ x\in \G_1\enspace\big|\enspace
dx+\frac{1}{2}[x,x] \in \OO_2\big\}\,.\]
It is clear that, for $x\in \G_1$ we have: $dx+\frac{1}{2}[x,x] \in
\OO_2$ if and only if $dx+\frac{1}{2}[x,x] \in
\OO$. This way the modular stack
$\M_1(\G,\OO):=\mc_1(\G,\OO)/G_0$ becomes a (locally-closed)
substack in
 $\M(\G,\OO)$.
\end{remark}

\subsection{}
Let $E=\bigoplus_{i\in \Z}\,E^i$ be a DG vector space
with differential $d: E^\hdot\to E^{\hdot+1}$,
and $(-,-): E\times E\map \C$
a $\C$-bilinear form such that for any homogeneous $x,y\in E$ we have
\beq{d_skew}
( dx,y) +(-1)^{\deg x}(x,dy)=0,
\eeq
and moreover $E^i\perp E^j$ whenever $i+j\neq 0.$

Put $E^*_i:=(E^{-i})^*,$ the dual of $E^{-i},$
and let $E^*:=\bigoplus_{i\in \Z}\,E^*_i.$
Dualizing the map $d$ makes $E^*$ a DG vector space.
The assignment $x\mapsto (x,-)$
gives rise to a morphism of DG vector spaces  $\kappa: E\to E^*$.
We say that the form   $(-,-)$
 is {\em non-degenerate} provided the  morphism $\kappa$ induces
an isomorphism $H^\hdot(\kappa): H^\hdot(E)\iso H^\hdot(E^*)$
on cohomology.

\subsection{}
Now, let $\G$ be a DGLA, as in \S\S\ref{sec1}-\ref{sec2},
and let $\be(-,-)$ be an {\it even} nondegenerate invariant 
bilinear form on $\G$, that is, a
$\C$-bilinear form $\G\times\G\to\C$ that restricts to a symmetric, 
resp. skew-symmetric, form on $\G_+$, resp.  on $\G_-,$
such that $\G_+\perp \G_-$, 
\eqref{d_skew} holds and 
for any homogeneous $x,y,z\in\G,$
one has:
\begin{equation}\l{be}
\be\bigl( [x,y],z\bigr) =
\be\bigl( x,[y,z]\bigr).
\end{equation}

\begin{theorem}\l{main}
Let $\G$ be a DGLA with an even nondegenerate invariant bilinear
form $\be$, and $\OO\sset \G_+$ an $\Ad G_+$-orbit.
Then, for any $x\in \M(\G,\OO)$, the form $\be$ induces a non-degenerate
2-form on $\DT_x\M(\G,\OO)$, the tangent space (at $x$) to the modular
stack.\footnote{See explanation in \S\ref{expl} below.}
 These 2-forms give rise to a symplectic structure on the stack
$\M(\G,\OO)$.
\end{theorem}
Thus,  a DGLA with an even nondegenerate invariant bilinear
form gives rise to symplectic  stacks.

\begin{remark} Morally, the main message of this paper is
that {\em all}  natural symplectic structures on 
 moduli spaces `arising in nature'
come from an appropriate
 (even) nondegenerate invariant bilinear
form on 
the DGLA that controls the moduli problem in question.
\end{remark}

We will now show that the construction
of the modular stack $\M(\G,\OO)$ is a special case 
of the Hamiltonian reduction. That will yield the proof of
 Theorem \ref{main}.

\subsection{Hamiltonian reduction construction.}\label{ham}
In the setup of Theorem \ref{main},
let $\om := \be|_{\G_{-}}$ and consider $\G_{-}$ as a
symplectic vector space equipped with the symplectic form $\om$.
Consider also the gauge action of the group $G_{+}$ on $\G_{-}$.
This is a symplectic action, and we claim that it has moment map
$$ \Phi: \G_{-} \to \G_{+}^{*} \stackrel{\be}{\simeq}\G_{+} : x 
\mapsto dx + \frac{1}{2}[x,x].$$

To this end, note that the differential of the map $\Phi$ at the point $x\in \G_{-}$
is given by
$$\Phi'_{x}:\ y\mto\Phi'_{x}(y) = dy+[x,y]\quad\forall\,y\in \G_{-}\cong T_x\G_-.$$ 
A direct computation yields
$ \Phi'_{x}(\xi_{a}(x)) = [a, \Phi(x)]\quad\forall\,a\in\G_{+}.$
It follows that $\Phi$ is $G_{+}$-equivariant. Next, by (\ref{be}),
we have 
\begin{equation}\label{def_mom}
\om(\xi_{a}(x), y) = \beta(a, \Phi'_{x}(y))
\quad\forall\,a\in\G_{+},\ x,y\in\G_{-}.
\end{equation}
It follows from the
 last equation combined  with the $G_+$-equivariance of $\Phi$
that the $G_+$-action is Hamiltonian with moment map $\Phi$.
Thus, $\M(\G,\OO)$ is precisely the Hamiltonian
reduction of $\G_{-}$ with respect to the orbit $\OO\subset\G_{+}\simeq\G_{+}^*$, see \cite{AM}.

\subsection{}\label{expl} In this subsection we explain the
meaning of the word `nondegenerate' in Theorem \ref{main}.
We will see that this
is closely related to the `self-dual nature' of the
Hamiltonian reduction construction.

In general, let $(X,\om)$ be a smooth symplectic variety
equipped with a Hamiltonian action of an algebraic
group $G$. Let $\g:=\Lie G$, fix a coadjoint orbit $\OO\sset\g^*$,
and write $\Phi: X \to \g^*$ for the moment map.

We consider the scheme $\Si:=\Phi\inv(\OO)$. If $\Phi'$,
the differential of $\Phi$, is surjective at any 
point $x\in X$ such that $\Phi(x)\in \OO$,
then $\Si$ is a smooth (locally closed) subscheme
of $X$. In general, if  $\Phi'$ is not necessarily
surjective, it is natural to view $\Si$ as a DG scheme,
which we denote $\Si\dg$.
The tangent space to this DG scheme at a closed point
$x\in \Si\dg$ is a DG vector space 
$\DT_x\Si\dg=T^0_x\Si\dg\bigoplus T^1_x\Si\dg,$ (concentrated in degrees $0$ and
$1$),
where
\begin{align}\label{DG1}
T^0_x\Si\dg:=T_xX,\quad &T^1_x\Si\dg:=\g^*,\\
&\text{with
differential}
\quad T^0_x\Si\dg=T_xX\stackrel{d}\too \g^*=T^1_x\Si\dg,\nonumber
\end{align}
given by the map $\Phi'_x: T_xX\to \g^*$.
We write $H^\hdot(\DT_x\Si\dg)$ for the cohomology
groups of the two-term complex above.
We have $H^0(\DT_x\Si\dg)=\Ker\Phi'_x$,
is the  Zariski tangent space to $\Si\dg$
(viewed as an ordinary scheme $\Si$),
and $H^1(\DT_x\Si\dg)=\g^*/\im\Phi'_x$,
is the space measuring the failure of
the differential of $\Phi$ to be surjective.

Next, we perform the Hamiltonian reduction 
and consider the quotient
$\Si/G$ (where $\Si$ is an ordinary scheme, not a DG scheme).
 If the $G$-action on  $\Si$ 
is free, then this quotient is a well-defined scheme again.
In general, for a not necessarily free action,
it is natural to consider
$\M:=\Si/G$ as a stack, cf. \cite{LMB},\cite{To}. Given $x\in \Si$,
write $\byy$ for the corresponding point in $\M$. Then,
the tangent space at $\byy$ to the  stack $\M=\Si/G$
is a DG vector space 
$\DT_\byy\M=T^{-1}_\byy\M\bigoplus T^0_\byy\M$
(concentrated in degrees $-1$ and
$0$), where
\begin{align}\label{DG2}
T^{-1}_\byy\M:=\g,\quad &T^0_\byy\M:=T_x\Si,\\
&\text{with
differential}
\quad T^{-1}_\byy\M=\g\stackrel{d}\too T_x\Si=T^0_\byy\M,\nonumber
\end{align}
given by the  derivative (at $1\in G$)
of the action-map $g\mto g(x)$.
Then, for the cohomology
groups, we have $H^0(\DT_\byy\M)=T_x\Si/\g\cdot x,$ is the normal
space to the $G$-orbit through $x\in\Si$,
and $H^{-1}(\DT_\byy\M)=\g^x$ is the Lie algebra
of the isotropy group of the point $x$, that measures
 the failure of the $\g$-action on $\Si$ to be infinitesimally-free.

Now, if we view (as has been explained earlier)
$\Si\dg:=\Phi\inv(\OO)$ as a DG scheme rather than an
ordinary scheme, then 
 the quotient
$\M\dg:=\Si\dg/G$ becomes a DG stack rather than an
ordinary stack.
To get the
tanget space of this DG stack,
we must combine formulas \eqref{DG1} and \eqref{DG2}
together. Thus, the tangent space to $\M\dg$
is  a DG vector space 
$\DT_\byy\M\dg=T^{-1}_\byy\M\dg\bigoplus T^0_\byy\M\dg\bigoplus T^1_\byy\M\dg,$
concentrated in degrees $-1,0,1$, 
such that the corresponding 3-term complex 
$T^{-1}_\byy\M\dg\to T^0_\byy\M\dg\to T^1_\byy\M\dg$ reads:
\beq{DG3}
\DT_\byy\M\dg: 
\xymatrix{\g\ar[rr]^<>(.5){\text{action}}&& T_xX
\ar[rr]^<>(.5){\Phi'_x}&&\g^*.}
\eeq

An important feature of \eqref{DG3} is that this complex is 
{\em self-dual}.
In more detail,   the definition of
moment map implies that the map $\Phi'_x: T_xX\map\g^*$
is obtained by composing the isomorphism $\om: T_xX\iso T^*_xX,$
induced by the symplectic form $\om,$ with the adjoint $a^\top:
T^*_xX\map \g^*$ of
the action-map
$a: \g\map T_xX$, that is,  one has: 
$$\Phi'_x=a^\top\ccirc\om:\ T_xX\iso T^*_xX\too  \g^*.$$
This shows that the DG vector space $(\DT_\byy\M\dg)^*$ is canonically
isomorphic to $\DT_\byy\M\dg$. The isomorphism
$\DT_\byy\M\dg\iso(\DT_\byy\M\dg)^*$ induces, of course, an isomorphism
of cohomology groups, and the non-degeneracy of symplectic 
form on the DG stack $\M\dg$ follows.

\subsection{Example: moduli of $G$-local systems.}
A typical (infinite-dimensional) example of the Theorem above
is the moduli space of bundles with flat connection
on a compact $C^\infty$-manifold $X$.

In more details, let $G$ be a  Lie group 
with Lie algebra $\g$, and
$P$ a principal $G$-bundle on $X$.
Let $\gp$ denote the associated vector bundle (with fiber $\g$)
corresponding to the adjoint representation $\Ad: G\to \GL(\g)$.
Thus, $\gp$ is a bundle of Lie algebras.

 Let  $\Omega^i(X, \gp)$ be 
 the vector space
of $C^\infty$-differential $i$-forms on $X$ with values in $\gp$.
The Lie bracket on $\gp$ combined with 
wedge-product of differential forms makes the graded space
$\G_P:=\bigoplus_i\,\Omega^i(X, \gp)$
a Lie super-algebra.

We are interested in the moduli space (or stack) $\M(P)$ of
 flat  $C^\infty$-connections on $P$, modulo gauge equivalence.
So, assume that $\M(P)$ is non-empty and choose
some  flat  $C^\infty$-connection $\nabla$ on $P$.
The connection induces a differential
$\nabla: \Om^\hdot (X, \gp)\map \Om^{\hdot+1} (X, \gp)$,
thus gives $\G_P$ the structure of a DGLA.

Any  connection on $P$ can be written in the form
$\nabla' = \nabla+\gamma$, for some
$\gamma\in\Omega^{1}(X,\gp)$. 
The  curvature of $\nabla'$
is $\nabla'\ccirc\nabla'=\nabla\gamma+\frac{1}{2}[\gamma,\gamma].$
Thus,  $\nabla'$
is flat if and only if 
$\gamma$ satisfies the Maurer-Cartan equation
$\nabla\gamma+\frac{1}{2}[\gamma,\gamma]=0$.
Thus, in the notation of Remark \ref{deg1}
(with $d:=\nabla$),
we have $\M(P)\cong \mc_1(\G_P,\{0\})/G_0\,
\bigl(=\M_1(\G_P,\{0\})\bigr),$
where $G_0:=\Aut(P)$ is the infinite-dimensional
group of gauge transformations.

To proceed further, we
assume $X$ to be compact oriented, and assume also that
there is a non-degenerate invariant symmetric bilinear form 
$\langle-,-\rangle: \g\times\g\to{\mathbb R}$.
The form $\langle-,-\rangle$
induces a nondegenerate pairing
$\langle-,-\rangle_{_P}: \gp\times\gp\map C^\infty(X)$.
It is straightforward to verify that the following
formula ($i=0,\ldots, n=\dim_{\mathbb{R}} X$):
$$\be:\ \xymatrix{\Omega^{i}(X,\gp)\times
\Omega^{n-i}(X,\gp)\ar[r]^<>(.5){\wedge}&
(\gp\otimes\gp)\otimes\Omega^{n}(X)\ar[rr]^<>(.5){\langle-,-\rangle_{_P}}&&
\Omega^{n}(X)\ar[r]^<>(.5){\int_X}&\mathbb{R}}
$$
gives a nondegenerate symmetric (even) 
bilinear form $\beta: \G_P\times \G_P\map \mathbb{R}.$
The form $\beta$ is  {\em invariant}, provided
the connection $\nabla$ was chosen so that
the paring $\langle-,-\rangle_{_P}$ is $\nabla$-horizontal,
i.e., such that $\langle \nabla  x,y\rangle_{_P}+\langle x,\nabla  y\rangle_{_P}=0,$
for any sections
$x,y\in \gp$, c.f.~\eqref{d_skew}.

Now if $\dim_{\mathbb{R}} X=2$, then
$\Omega^{i}(X,\gp)=0$ for $i>2$, and therefore
we have $\mc(\G_P,\{0\})=\mc_1(\G_P,\{0\})$.
Hence the  construction of \S\ref{ham}
applied to the DGLA $\G=\G_P$ gives a symplectic structure
on the stack $\M(P)=\mc_1(\G_P,\{0\})/G_0=\mc(\G_P,\{0\})/G=\M(\G_P,\{0\})$.

For $\dim_{\mathbb{R}} X>2,$
our construction only gives a  symplectic structure
on $\M(\G_P,\{0\})\supsetneq \M(P)$.
However, it is known (see e.g., \cite{Kar})
that if the  Hard Lefschetz theorem holds for $X$,
then $\M(P)$ is in effect a  symplectic  substack in $\M(\G_P,\{0\})$.

\begin{remark} For a trivial $G$-bundle $P$, the space $\M(P)$
may be identified with $\Hom(\pi_1(X),G)/\Ad G$,
the moduli space of $G$-local systems on $X$. The symplectic
structure on this space has been studied by many authors,
see e.g. \cite{Go},\cite{We}, \cite{Kar}.
\end{remark}

\section{$L_\infty$-algebra version}
\subsection{}
Let $\G=\G_+\oplus\G_-$ be a $\Z/2\Z$-graded $L_\infty$-algebra, 
see e.g. [LM], [MSS]. 
Let $\be$ be an even nondegenerate invariant bilinear 
form on $\G$. Thus, for any homogeneous $x_{1},\ldots,x_{n+1} \in\G$,
one has:
$$ \be([x_{1},\ldots,x_{n}],x_{n+1}) = (-1)^{n(\deg x_{1}+1)}
\be(x_{1}, [x_{2}, \ldots,x_{n+1}]).$$
Define $$ \Phi:\G_{-}\to\G_{+}:x\mapsto dx+\frac{1}{2!}[x,x]
+\frac{1}{3!}[x,x,x]+\cdots.$$ 
The differential of the map $\Phi$ at the point $x\in \G_{-}$
is given by
$$\Phi'_{x}:\ y\mto\Phi'_{x}(y) = dy+[x,y]+\frac{1}{2!}[x,x,y]+\cdots,\quad
\forall y\in\G_{-}. $$
For any $a\in\G_{+}$ and $x\in\G_{-}$, let
\beq{ser1} \xi_{a}(x) :=
-\bigl(da+[x,a]+\frac{1}{2!}[x,x,a]+\cdots\bigr).
\eeq
Observe that we have
\begin{equation}\l{mmm}
\beta(\xi_{a}(x), y) = \beta(a, \Phi'_{x}(y))
\quad\forall\,a\in\G_{+}\,,x,y\in\G_{-}.
\end{equation}
By [La, Appendix B], we also have
\beq{ser2}
\Phi'_{x}(\xi_{a}(x)) = [a,\Phi(x)]+[a,\Phi(x),x]+
\frac{1}{2!}[a, \Phi(x), x,x]+\cdots\quad
\forall\,a\in\G_{+}\,,x\in\G_{-}.
\end{equation}

\noindent
\pb{We say that two elements $x_{0},x_{1} \in \G_{-}$ are \emph{gauge equivalent}
if there exists a path $a(t)\in\G_{+}$ and
a path $x(t)\in\G_{-}$ such that
$ x'(t) = \xi_{a(t)}(x(t))$ and $x(0)=x_{0},\,x(1)=~x_{1}$.}

\noindent
\pb{We say that two elements $b_{0}, b_{1}\in\G_{+}$ are \emph{adjoint equivalent}
if there exists two paths $a(t), b(t)\in\G_{+}$ and
a path $x(t)\in\G_{-}$ such that $b(0)=b_{0}, b(1)=b_{1}$ and}
\begin{equation}\l{adjoint}
b'(t)=[a(t),b(t)]+[a(t),b(t),x(t)]+\frac{1}{2!}
[a(t),b(t),x(t),x(t)] +\cdots. 
\end{equation}

\begin{remark}
If $\G=\bigoplus_{i\in\Z}\,\G_{i}$ is a $\Z$-graded $L_{\infty}$-algebra,
then we say that $b_{0}, b_{1}\in\G_{2}$ are adjoint equivalent
if there exists a path $a(t)\in\G_{0}$, a path
$b(t)\in\G_{2}$ and a path $x(t)\in\G_{1}$ such that $b(0)=b_{0},
b(1)=b_{1}$ and (\ref{adjoint}) is satisfied.
\end{remark}

\subsection{}
If $\OO\sset\G_+$ is an adjoint equivalence class, then we define
$$\mc(\G,\OO):= \{x\in\G_-\,|\,\Phi(x)\in\OO\},\quad\M(\G,\OO):=\mc(\G,\OO)/
\textrm{gauge equivalence}.$$

\begin{remark} This definition of the objects $\mc(\G,\OO)$ and
$\M(\G,\OO)$ involves solving differential equations, as well
as formulas \eqref{ser1}-\eqref{ser2} which
contain infinite
series. One way to make sense of the definition above is
to use a topological setting and to prove the convergence
of all the 
series  that arise. 

If $\OO=\{0\}$, then there is an alternative, purely algebraic, approach
based on the language of {\em formal schemes}.
Specifically, let $\wh{\G}_-$ be the
 formal 
completion of the vector space $\G_-$ at the origin.
The coordinate ring of $\wh{\G}_-$ is the ring
$\C[[\G_-]]$ of formal power series, so that the series for $\Phi$ gives a well-defined
morphism $\wh{\G}_-\to\G_+$.
Thus, $\Phi\inv(0)\sset \wh{\G}_-$ is a well-defined closed subscheme.
Further, gauge equivalence gives a well-defined
pro-algebraic groupoid acting on this subscheme,
and one puts $\wh{\M}(\G,\{0\}):=\Phi\inv(0)/\text{\it gauge equivalence}$,
a pro-algebraic stack.
\end{remark}

One can extend the concept of Hamiltonian reduction
from Hamiltonian group-actions to Hamiltonian
{\em groupoid-actions}. This way, using formula
(\ref{mmm}), one derives the following

\begin{theorem}\l{main_infty}
Let $\G=\G_+\oplus\G_-$ be an $L_\infty$-algebra
with an invariant nondegenerate even bilinear
form $\be$, and $\OO\sset \G_+$ an adjoint equivalence class.
Then, for any $x\in \M(\G,\OO)$,
the form $\be$ induces a non-degenerate 2-form
on $\DT_x\M(\G,\OO)$, the tangent space (at $x$) to the modular stack.
These 2-forms give rise to a symplectic structure on
$\M(\G,\OO)$.
\end{theorem}

\section{From moment map to Maurer-Cartan equations}
\subsection{}
Let $(V,\om)$ be a finite dimensional symplectic vector
space, and $\C[V]$ the  algebra of polynomial
functions on $V$ viewed as a
Poisson algebra with Poisson bracket $\{-,-\}$
corresponding to the symplectic structure.

Fix  a
 finite dimensional Lie algebra $\g$. Let
 $H: \g\map \C^{>0}[V]=\bigoplus_{i\geq 1} \C^i[V]
\,,\, a\mapsto H_a,$ be a Lie algebra homomorphism.
We get, tautologically, 
a (non-linear) Hamiltonian $\g$-action on the vector space
 $V$ with polynomial moment map
$\Phi: V\to \g^*\,,\,\Phi(v): a \mapsto H_a(v).$

\subsection{}
Introduce a $\Z$-graded vector space
$$\G= \g\oplus V \oplus \g^*,\qquad
\G_0:=\g\,,\,\G_1:=V\,,\,\G_2:=\g^*,$$
and write $\G_+:=\G_0\oplus\G_2$ and $\G_-:=\G_1$.
The canonical pairing $\g\times\g^*\to\C$ gives
a non-degenerate symmetric bilinear form on $\G_+$.
Combined with the symplectic form $\om$ on $V$ this
gives an even non-degenerate bilinear form $\be$ on
the super-space $\G=\G_+\oplus\G_-$. 

Further, let $\Phi=\Phi_1 + \Phi_2+\ldots\, :
V \too \g^*$ be  an expansion of the moment
map $\Phi$  into homogeneous components
(by assumption there is no constant term).
For each $i\geq 1$, the map $\Phi_i$
gives rise, via the canonical isomorphism
$\C^i[V]\simeq (\sym^iV)^*,$ to a {\em linear} map 
$\tilde\Phi_i: \sym^iV \to \g^*,$ such that
$\tilde\Phi_i(y^i)=\Phi_i(y)\,,\,\forall y\in V$.

We now define the following maps:

\noindent
\pb{$\G_0 \otimes \G\map \G$ given, for any $u=x\oplus v\oplus\lambda\in
\g\oplus V\oplus \g^*,$ by $a\otimes u\mapsto [a,u]:=\ad a(x)\oplus
a(v)\oplus\ad^*a(\lambda);$}

\noindent
\pb{$\sym^i(\G_-)\to \G_2\sset \G_+\,,\,y_1 y_2\ldots y_i\mapsto
[y_1,y_2,\ldots,y_i]:=\tilde\Phi_i(y_1 y_2\ldots y_i),
$ for each $i\geq 2;$}

Further, let $d:=\Phi_1: \ \G_-\to \G_+$, and
define $d: \G_+ \to\G_-$ to be the zero-map.
\begin{proposition} The above defined maps give $\G$ 
 an $L_\infty$-structure (with all brackets that were not specified
above being set equal to zero), such that
$\be$ becomes an invariant form.
\end{proposition}
\begin{proof} Straightforward computation.
\end{proof}
\vskip -5mm
We observe that two elements in $ \G_2=\g^*$
are adjoint equivalent (in the sense of \S2) if and only if they belong to the same
coadjoint orbit of the group $G=G_0$ (corresponding to the Lie algebra
$\g$)
acting on $\g^*$. Given such an
orbit $\OO\sset \G_2
=\g^*$,
we see that $\mc(\G,\OO)=\Phi^{-1}(\OO)$.
Hence, for the Maurer-Cartan stack we get $\M(\G,\OO)=
\Phi^{-1}(\OO)/\Ad G_0,$ is the standard Hamiltonian reduction
of $V$ over $\OO$.

\noindent
{\bf Remarks.}$\enspace$\vi
Note that the construction above may be thought of as a {\em non-linear}
version of the construction of the tangent space as a DG vector space,
given in \S\ref{ham}.

\vii In the special case 
 $H: \g\map \C^2[V]$ ({\em quadratic Hamiltonians}),
we have $\Phi=\Phi_2$, and the  $L_\infty$-structure above
reduces to an ordinary Lie super-algebra structure
on $\G=\g\oplus V\oplus\g^*.$
The {\it symmetric} bracket $\G_-\times\G_- \to \G_+$
is in this case provided by the map
$\Phi$ viewed as a linear map: $\sym^2(\G_-)=\sym^2(V)\to\g^*=\G_2.$
\footnote{We are grateful to E. Getzler for pointing 
out a similarity between our construction
and the results of  Goldman and Millson, c.f. \cite[p. 499]{GM2}.}
This special case was also implicit in \cite{Ko}. Notice
that a Lie algebra homomorphism $\g \map \sym^2 V \simeq
\mathfrak{sp}(V) $
has a natural extension to a Lie algebra homomorphism 
$\nu:\G_+=\g\oplus\g^* \map \mathfrak{sp}(V)$
 by sending $\g^*$ to zero.
We observe that $\nu$ is of super Lie type in Kostant's 
terminology, i.e. the condition on the Casimir in 
\cite[Theorem 0.1]{Ko} holds (trivially) for this $\nu$.

\viii J. Stasheff informed us that it is possible to 
define Hamiltonian reductions with respect to
 (infinitesimal) actions of an $L_\infty$-algebra, and to
extend the construction of \S3 to such an $L_\infty$-setup.
He also pointed out to us the relationship of our construction
to the classical BRST complex. Namely, an $L_\infty$-structure
on $\G$ is equivalent to a square zero derivation 
on the free super-commutative algebra
$S^{>0}(\G[1]^*)$ generated by $\G[1]^*$, 
see e.g. [LM]. Here, $\G[1]$ is the super 
vector space  with $\G[1]_+ = \G_-$ and $\G[1]_- = \G_+$.
For the $L_\infty$-algebra $\G$ in Proposition 3.2.1, we thus
obtained a differential on $S^{>0}(\G[1]^*)$ which turns out to
be the classical BRST operator on $\bigwedge\g^* \otimes \bigwedge\g
\otimes \sym V^*$ defined in \cite[p. 57]{KS}. 
To see this, it suffices to
note that the generators of $S^{>0}(\G[1]^*)$ are
$$ \G[1]^*_- = \g^* \oplus \g \ , \quad \quad 
\G[1]^*_+ = V^*,$$
and the differential is defined on these generators by taking the
sum of the maps 
$$H: \g \map \sym^{>0}V^*  \quad \textrm{and} \quad
J: \G[1]^* \map \g^* \otimes \G[1]^*, $$
where $J$ is obtained from  dualizing the map $\G_0 \otimes \G 
\map \G: a\otimes u \mapsto [a,u]$.

\setcounter{equation}{0}
\footnotesize{

Department of Mathematics, Massachusetts Institute of Technology,
Cambridge, MA 02139, USA;\\
\hphantom{x}\quad\, {\tt wlgan@math.mit.edu}

\smallskip

Department of Mathematics, University of Chicago, 
Chicago, IL 60637, USA;\\ 
\hphantom{x}\quad\, {\tt ginzburg@math.uchicago.edu}}

\end{document}